\documentclass[12pt]{article}
\usepackage{amsmath}
\usepackage{amssymb}
\usepackage{vatola}

\def\q{\quad}
\def\qq{\qquad}
\def\qtq#1{\q\t{#1}\q}
\def\mod#1{\ (\text{\rm mod}\ #1)}
\def\t{\text}
\def\f{\frac}
\def\e{\equiv}
\def\b{\binom}

\def\ap{\langle a\rangle_p}
\def\sls#1#2{(\f{#1}{#2})}
 \def\ls#1#2{\big(\f{#1}{#2}\big)}
\def\Ls#1#2{\Big(\f{#1}{#2}\Big)}

\let \pro=\proclaim
\let \endpro=\endproclaim

\begin{document}
\par\q\par\q
 \centerline {\bf
Congruences involving Franel and Catalan-Larcombe-French numbers}
$$\q$$
\centerline{Zhi-Hong Sun} $$\q$$ \centerline{School of Mathematical
Sciences, Huaiyin Normal University,} \centerline{Huaian, Jiangsu
223001, P.R. China} \centerline{Email: zhihongsun@yahoo.com}
\centerline{Homepage: http://www.hytc.edu.cn/xsjl/szh}
 \abstract{Let $\{f_n\}$ be the Franel numbers given by
  $f_n=\sum_{k=0}^n\binom  nk^3$, and let $p>5$ be a prime. In this paper
  we mainly determine  $\sum_{k=0}^{p-1}
 \binom{2k}k\frac{f_k}{m^k}\pmod p$ for $m=5,-16,16,32,-49,50,96$.
 Let $S_n=\sum_{k=0}^n\binom nk\binom{2k}k\binom{2n-2k}{n-k}$.
  We also determine
$\sum_{k=0}^{p-1}
 \binom{2k}k\frac{S_k}{m^k}\pmod p$ for $m=7,16,25,32,64,160,800,1600,
 156832$.
 \par\q
 }
 \endabstract
\let\thefootnote\relax \footnotetext {The author is supported by
the Natural Science Foundation of China (grant No. 11371163).}

\section*{1. Introduction}
\par\q  Let $[x]$ be the greatest integer not exceeding $x$,
and let $\sls ap$ be the Legendre symbol. For a prime $p$ let $\Bbb
Z_p$ be the set of rational numbers whose denominator is not
divisible by $p$. For positive integers $a,b$ and $n$, if
$n=ax^2+by^2$ for some integers $x$ and $y$, we briefly write that
$n=ax^2+by^2$.
\par In 1894 J. Franel [F] introduced the following Franel numbers
$\{f_n\}$:
$$f_n=\sum_{k=0}^n\b nk^3\q(n=0,1,2,\ldots).$$
The first few Franel numbers are as below:
$$f_0=1,\ f_1=2,\ f_3=10,\ f_4=56,\
f_5=346,\ f_6=2252,
 \ f_7=15184.$$
 It is known that
 $$(n+1)^2f_{n+1}=(7n^2+7n+2)f_n+8n^2f_{n-1}\ (n\ge 1).$$
 \par Let $p$ be an odd prime and
 $m\in\Bbb Z$ with $m\not\e 0\mod p$. In [S6],
  the author made many conjectures on $\sum_{k=0}^{p-1}
 \b{2k}k\f{f_k}{m^k}\mod {p^2}$. For example, for any odd prime $p$,

$$\sum_{k=0}^{p-1}\b{2k}k\f{f_k}{(-16)^k}\e \cases 4x^2-2p\mod {p^2}&\t{if $p=x^2+y^2\e 1\mod{12}$ with $6\mid y$,}
\\2p-4x^2\mod{p^2}&\t{if $p=x^2+y^2\e 1\mod{12}$ with $6\mid x-3$,}
\\4\sls {xy}3xy\mod{p^2}&\t{if $p=x^2+y^2\e 5\mod{12}$,}
\\0\mod{p^2}&\t{if $p\e 3\mod 4$.}\endcases$$
 In [Gu], J.W. Guo proved that
 $$\sum_{k=0}^{p-1}\b{2k}k\f{f_k}{(-16)^k}\e 0\mod p\q\t{for}\q p\e
 3\mod 4$$ and
 $$\sum_{k=0}^{p-1}\f{3k+1}{(-16)^k}\b{2k}kf_k\e
 (-1)^{\f{p-1}2}p\mod {p^3},$$
 where the second congruence modulo $p^2$
 was conjectured by the author in [S6]. We note that $p\mid \b{2k}k$
 for $k=\f{p+1}2,\ldots,p-1$.
 In [Su4, Su5], the author's
 brother Z.W. Sun investigated congruences for Franel numbers. In
 particular, he showed that for any odd prime $p$,
 $$\sum_{k=0}^{p-1}\b{2k}k\f{f_k}{(-4)^k}\e
 \sum_{k=0}^{p-1}\f{\b{2k}k^3}{16^k}\mod{p^2}.$$
 By [S3, Theorems 3.3 and 3.4],
$$\sum_{k=0}^{p-1}\f{\b{2k}k^3}{16^k}\e\cases
4x^2-2p\mod{p^2}&\t{if $p\e 1\mod 3$ and so $p=x^2+3y^2$,}
\\0\mod{p^2}&\t{if $p\e 2\mod 3$.}
\endcases$$ Thus,
$$\sum_{k=0}^{p-1}\b{2k}k\f{f_k}{(-4)^k}\e\cases
4x^2-2p\mod{p^2}&\t{if $p\e 1\mod 3$ and so $p=x^2+3y^2$,}
\\0\mod{p^2}&\t{if $p\e 2\mod 3$.}
\endcases$$
\par For any nonnegative integer n let
$$\aligned &A_n=\sum_{k=0}^n\b nk^2\b{n+k}k^2,\q D_n=\sum_{k=0}^n\b{2k}k\b{2n-2k}{n-k}\b nk^2,
\\& a_n=\sum_{k=0}^n\b nk^2\b{2k}k,\q
b_n=\sum_{k=0}^{[n/3]}\b{2k}k\b{3k}k\b n{3k}\b{n+k}k(-3)^{n-3k},
\\&S_n=\f{P_n}{2^n}=\sum_{k=0}^n\b nk\b{2k}k\b{2n-2k}{n-k}
=\sum_{k=0}^{[n/2]}\b n{2k}\b{2k}k^24^{n-2k}.
\endaligned\tag 1.1$$
Here $\{A_n\}$ is called Ap\'ery numbers since Ap\'ery [Ap] used it
to prove $\zeta(3)$ is irrational in 1979, $\{D_n\}$ is called Domb
numbers, $\{b_n\}$ is called Almkvist-Zudilin numbers, and $\{P_n\}$
is called Catalan-Larcombe-French numbers. See [CCL], [CV], [Co],
[D], [JV], [Su6] and [Z].  Such sequences appear as coefficients in
various series for $1/\pi$, for example,
$$\sum_{k=0}^{\infty}\f{9k+2}{50^k}\b{2k}kf_k=\f{25}{2\pi},
\q \sum_{k=0}^{\infty}\f{5k+1}{64^k}D_k=\f 8{\sqrt 3\pi},\q
\sum_{n=0}^{\infty}\f{4k+1}{81^k}b_k=\f{3\sqrt 3}{2\pi}. $$
\par Let $p>3$ be a prime, $m\in\Bbb Z_p$ and
$m\not\e 0,-4,-8\mod p$. In this paper, we show that
$$\sum_{k=0}^{p-1}\b{2k}k\Ls m{(m+8)^2}^kf_k\e
\sum_{k=0}^{p-1}\b{2k}k^2\b{3k}k\Ls m{(m+4)^3}^k\mod p.$$ Let
$x\in\Bbb Z_p$, $x\not\e 0,-1,-\f 13\mod p$ and $\ls
{9x^2+14x+9}p=1$. We also show that
$$\sum_{k=0}^{p-1}\b{2k}k\Ls x{9x^2+14x+9}^kf_k\e
\sum_{k=0}^{p-1}\b{2k}k^2\b{4k}{2k}\Ls x{9(1+3x)^4}^k\mod p.$$
 As consequences we determine $\sum_{k=0}^{p-1}
 \b{2k}k\f{f_k}{n^k}\mod p$ for $n=5,-16,16,32,-49,50,96$ and $\sum_{k=0}^{p-1}
 \b{2k}k\f{a_k}{4^k}\mod p$. As examples, for any prime $p>5$ we have
$$\align &\sum_{k=0}^{p-1}\b{2k}k\f{f_k}{(-16)^k}\e 4x^2\mod p
\qtq{for}p=x^2+9y^2, \\&\sum_{k=0}^{p-1}\b{2k}k\f{f_k}{16^k}\e
4x^2\mod p\qtq{for}p=x^2+5y^2,
\\&\sum_{k=0}^{p-1}\b{2k}k\f{f_k}{5^k}\e 4x^2\mod
p\qtq{for}p=x^2+15y^2, \\&\sum_{k=0}^{p-1}\b{2k}k\f{f_k}{32^k}\e
4x^2\mod p\qtq{for}p=x^2+6y^2.\endalign$$
 Thus we partially solve some conjectures in [S6].
\par In [Su6] Z.W. Sun introduced
$$S_n(x)=\sum_{k=0}^n\b nk\b{2k}k\b{2n-2k}{n-k}x^k
\q (n=0,1,2,\ldots)$$ and used it to establish new series for
$1/\pi$. Note that $S_n(1)=S_n$ is essentially the
Catalan-Larcombe-French number. In [JV], Jarvis and Verrill gave
some congruences for $P_n=2^nS_n$. In Section 3 we establish some
new identities involving $S_n$. For example,
$$\sum_{k=0}^n\b nk(-1)^k\f{S_k}{8^k}=\f{S_n}{8^n}
\qtq{and}
\sum_{k=0}^{2n}\b{2n}k\b{2n+k}k(-8)^{2n-k}S_k=(-1)^n\b{2n}n^3.$$
  Let $p$ be an odd
prime, $n\in\Bbb Z_p$ and $n\not\e 0,-16\mod p$. In Section 3 we
also prove that
$$\sum_{k=0}^{p-1}\b{2k}k\f{S_k}{(n+16)^k}\e \Ls {n(n+16)}p\sum_{k=0}^{p-1}
\f{\b{2k}k^2\b{4k}{2k}}{n^{2k}}\mod p.$$
 As consequences we determine $\sum_{k=0}^{p-1}\b{2k}k\f{S_k}{m^k}\mod
p$ for $m=7,16,25,32,64,160,800,1600,$ $156832$. For example, for
any prime $p>7$,
$$\sum_{k=0}^{p-1}\b{2k}k\f{S_k}{7^k}
\e\cases 4x^2\mod p&\t{if $p\e 1,2,4\mod 7$ and so $p=x^2+7y^2$,}
\\0\mod p&\t{if $p\e 3,5,6\mod 7$.}
\endcases$$
\par Let $p$ be an odd prime, $n,x\in\Bbb Z_p$ and
$n(n+4x)\not\e 0\mod p$. In Section 4 we show that
$$\sum_{k=0}^{p-1}\b{2k}k\f{C_k(x)}{(n+4x)^k}
\e\Ls {n(n+4x)}p\sum_{k=0}^{p-1}\f{\b{2k}k\b{3k}k\b{6k}{3k}}
{n^{3k}}\mod p,$$ where
$$C_n(x)=\sum_{k=0}^{[n/3]}\b{2k}k\b{3k}k\b n{3k}x^{n-3k}.$$
\par In this paper we also pose some conjectures for congruences
involving $f_n$ or $S_n$. See Conjectures 2.1-2.2 and Conjectures
3.1-3.4.
\section*{2.  Congruences involving $\{f_n\}$}
 \pro{Lemma 2.1} Let $p$ be an odd prime, $u\in\Bbb Z_p$ and
 $u\not\e  1\mod p$. For any $p$-adic sequences $\{c_k\}$ we
 have
 $$\sum_{k=0}^{p-1}\b{2k}k\Ls u{(1-u)^2}^kc_k
 \e\sum_{n=0}^{p-1}u^n\sum_{k=0}^n\b nk\b {n+k}kc_k\mod p.$$
 \endpro
 Proof. Note that $\b{-x}k=(-1)^k\b{x+k-1}k$ and $\b nk\b{n+k}k
 =\b{2k}k\b{n+k}{2k}$. Using Fermat's little theorem we deduce that
$$\align
&\sum_{k=0}^{p-1}\b{2k}k\Ls u{(1-u)^2}^kc_k \\&\e
\sum_{k=0}^{(p-1)/2}\b{2k}kc_ku^k(1-u)^{p-1-2k}
=\sum_{k=0}^{(p-1)/2}\b{2k}kc_ku^k\sum_{r=0}^{p-1-2k}\b{p-1-2k}r(-u)^r
\\&=\sum_{n=0}^{p-1}u^n\sum_{k=0}^n\b
{2k}kc_k(-1)^{n-k}\b{p-1-2k}{n-k} =\sum_{n=0}^{p-1}u^n\sum_{k=0}^n\b
{2k}kc_k\b{n+k-p}{n-k}
\\&\e \sum_{n=0}^{p-1}u^n\sum_{k=0}^n\b
{2k}kc_k\b{n+k}{n-k} =\sum_{n=0}^{p-1}u^n\sum_{k=0}^n\b nk\b
{n+k}kc_k\mod p.\endalign$$ Thus the lemma is proved.

\pro{Lemma 2.2} Let $p>3$ be a prime and
$c_0,c_1,\ldots,c_{p-1}\in\Bbb Z_p$. Then
$$\sum_{n=0}^{p-1}\sum_{k=0}^n\b nk\b{n+k}kc_k
\e \sum_{k=0}^{p-1}\f{p}{2k+1}(-1)^kc_k\mod{p^3}.$$
\endpro
Proof. Since
$$\align \sum_{k=0}^m\b xk(-1)^k&=\sum_{k=0}^m\b{x-1}k(-1)^k+
\sum_{k=1}^m\b{x-1}{k-1}(-1)^k
\\&=\sum_{r=0}^m\b{x-1}r(-1)^r-\sum_{r=0}^{m-1}\b{x-1}r(-1)^r
\\&=\b{x-1}m(-1)^m=\b{m-x}m,\endalign$$
and $\b nk\b{n+k}k=\b{2k}k\b{n+k}{2k}$ we see that
$$\align &\sum_{n=0}^{p-1}\sum_{k=0}^n\b nk\b{n+k}kc_k
\\&=\sum_{k=0}^{p-1}\sum_{n=k}^{p-1}\b{2k}k\b{n+k}{2k}c_k
=\sum_{k=0}^{p-1}\b{2k}kc_k\sum_{r=0}^{p-1-k}\b{2k+r}{2k}
\\&=\sum_{k=0}^{p-1}\b{2k}kc_k\sum_{r=0}^{p-1-k}\b{-2k-1}r(-1)^r
=\sum_{k=0}^{p-1}\b{2k}kc_k\b{p+k}{p-1-k}
\\&=\sum_{k=0}^{p-1}\b{2k}kc_k\b{p+k}{2k+1}
=\sum_{k=0}^{p-1}\f p{2k+1}c_k\f{(p^2-1^2)\cdots(p^2-k^2)}{k!^2}.
\endalign$$ Thus,
$$\aligned &\sum_{n=0}^{p-1}\sum_{k=0}^n\b nk\b{n+k}kc_k
\\&=\f{(p^2-1^2)\cdots(p^2-\sls{p-1}2^2)}{\sls{p-1}2!^2}c_{\f{p-1}2}
 +p\sum\Sb
k=0\\k\not=(p-1)/2\endSb^{p-1}\f{c_k}{2k+1}\cdot
\f{(p^2-1^2)\cdots(p^2-k^2)}{k!^2}
\\&\e (-1)^{\f{p-1}2}\Big(1-p^2\sum_{r=1}^{(p-1)/2}\f
1{r^2}\Big)c_{\f{p-1}2}+p\sum\Sb
k=0\\k\not=(p-1)/2\endSb^{p-1}\f{(-1)^kc_k}{2k+1}\mod{p^3}.
\endaligned$$
It is known that $\sum_{r=1}^{(p-1)/2}\f 1{r^2}\e 0\mod p$. Thus the
result follows.
\par{\bf Example 2.1.} Let $\{P_n(x)\}$ be the famous Legendre
polynomials. Then Murphy proved that
$$P_n(x)=\sum_{k=0}^n\b nk\b{n+k}k\Ls{x-1}2^k.$$
Thus, applying Lemma 2.2 we see that for any prime $p>3$ and
$x\in\Bbb Z_p$,
$$\sum_{n=0}^{p-1}P_n(x)\e \sum_{k=0}^{p-1}\f p{2k+1}\Ls{1-x}2^k
\mod {p^3}.\tag 2.1$$

\pro{Lemma 2.3 ([CTYZ, (2.19), p.1305 and (2.27)]} Let $n$ be a
nonnegative integer. Then
$$A_n=\sum_{k=0}^n\b nk\b{n+k}kf_k=\sum_{k=0}^n\b
nk\b{n+k}k(-1)^{n-k}a_k$$ and
$$\f{D_n}{8^n}=\sum_{k=0}^n\b nk\b{n+k}k\f{f_k}{(-8)^k}.$$
\endpro
\par Lemma 2.3 can be verified straightforward by using Maple and the method in [CHM]
to compare the recurrence relations for both sides.

\pro{Theorem 2.1} Let $p$ be an odd prime, $m\in\Bbb Z_p$ and
$m\not\e 0,-4,-8\mod p$. Then
$$\sum_{k=0}^{p-1}\b{2k}k\Ls m{(m+8)^2}^kf_k\e
\sum_{k=0}^{p-1}\b{2k}k^2\b{3k}k\Ls m{(m+4)^3}^k\mod p.$$
\endpro

 Proof. Taking
$c_k=\f{f_k}{(-8)^k}$ in Lemma 2.1 and then applying Lemma 2.3 we
see that for $u\in\Bbb Z_p$ with $u\not\e 1\mod p$,
$$\sum_{k=0}^{p-1}\b{2k}k\Ls u{(1-u)^2}^k\f{f_k}{(-8)^k}
\e\sum_{n=0}^{p-1}u^n\f{D_n}{8^n}\mod p.$$ Now substituting $u$ with
$-\f 8m$ in the above formula we deduce that
$$\sum_{k=0}^{p-1}\b{2k}k\Ls m{(m+8)^2}^kf_k\e
\sum_{n=0}^{p-1}\f{D_n}{(-m)^n}\mod p.\tag 2.2$$ By [S8, Theorem
3.1],
$$\sum_{n=0}^{p-1}\f{D_n}{(-m)^n}\e
\sum_{k=0}^{p-1}\b{2k}k^2\b{3k}k\Ls m{(m+4)^3}^k\mod p.$$ Thus the
theorem is proved.

\pro{Theorem 2.2} Let $p>5$ be a prime. Then
$$\sum_{k=0}^{p-1}\b{2k}k\f{f_k}{50^k}
\e\cases 4x^2\mod p&\t{if $p\e 1\mod 3$ and so $p=x^2+3y^2$,}
\\0\mod p&\t{if $p\e 2\mod 3$.}
\endcases$$
\endpro
Proof. Taking $m=2$ in Theorem 2.1 we see that
$$\sum_{k=0}^{p-1}\b{2k}k\f{f_k}{50^k}
\e \sum_{k=0}^{p-1}\f{\b{2k}k^2\b{3k}k}{108^k}\mod p.$$ From [M] and
[Su2] we know that
$$\sum_{k=0}^{p-1}\f{\b{2k}k^2\b{3k}k}{108^k}
\e \cases 4x^2-2p\mod {p^2}&\t{if $p\e 1\mod 3$ and so
$p=x^2+3y^2$,}
\\0\mod {p^2}&\t{if $p\e 2\mod 3$.}
\endcases$$
Thus the result follows.

\pro{Theorem 2.3} Let $p$ be a prime with $p\e \pm 1\mod 8$. Then
$$\sum_{k=0}^{p-1}\b{2k}k\f{f_k}{32^k}
\e\cases 4x^2\mod p&\t{if $p\e 1,7\mod {24}$ and so $p=x^2+6y^2$,}
\\0\mod p&\t{if $p\e 17,23\mod{24}$.}
\endcases$$
\endpro
Proof. Taking $m=8$ in Theorem 2.1 we see that
$$\sum_{k=0}^{p-1}\f{\b{2k}kf_k}{32^k}\e
\sum_{n=0}^{p-1}\f{\b{2k}k^2\b{3k}k}{216^k}\mod p.$$ Now applying
[S2, Theorem 4.5] we deduce the result.

\pro{Theorem 2.4} Let $p$ be a prime with $p\e \pm 1\mod 5$. Then
$$\align
\sum_{k=0}^{p-1}\b{2k}k\f{f_k}{(-49)^k} \e\cases 4x^2\mod p&\t{if
$p\e 1,19\mod {30}$ and so $p=x^2+15y^2$,}
\\0\mod p&\t{if $p\e 11,29\mod{30}$.}
\endcases\endalign$$
\endpro
Proof. Taking $m=-1$ in Theorem 2.1 we see that
$$\sum_{k=0}^{p-1}\b{2k}k\f{f_k}{(-49)^k}\e
\sum_{k=0}^{p-1}\f{\b{2k}k^2\b{3k}k}{(-27)^k}\mod p.$$ Now applying
[S2, Theorem 4.6] we deduce the result.
 \pro{Theorem 2.5} Let $p$ be
a prime such that $p\e 1,19\mod{30}$ and so $p=x^2+15y^2$. Then
$$\sum_{k=0}^{p-1}\b{2k}k\f{f_k}{5^k}
\e 4x^2\mod p.$$
\endpro
Proof. Let $t\in\{1,2,\ldots,\f{p-1}2\}$ be given by $t^2\e -15\mod
p$ and $m=(-11+3t)/2$. Then $\f{64}m\e \f{-11-3t}2\mod p$ and so
$$\f{(m+8)^2}m=16+m+\f{64}m\e 16+\f{-11+3t}2+\f{-11-3t}2=5\mod p.$$
We also have
$$\f{(m+4)^3}m\e \f{\sls{-3+3t}2^3}{\f{-11+3t}2}\e -27\mod p.$$
Now applying Theorem 2.1 and [S2, Theorem 4.6] we deduce that
$$\sum_{k=0}^{p-1}\b{2k}k\f{f_k}{5^k}\e
\sum_{k=0}^{p-1}\f{\b{2k}k^2\b{3k}k}{(-27)^k}\e 4x^2\mod p.$$ This
proves the theorem.

\par{\bf Remark 2.1} Let $p$ be an odd prime. Taking $m=-16$ in
Theorem 2.1 we deduce the congruence for $\sum_{k=0}^{p-1}
\b{2k}k\f{f_k}{(-4)^k}\mod p$.

 \pro{Theorem 2.6} Let p be an odd prime  and $u\in\Bbb Z_p$.
\par $(\t{\rm i})$ If $u\not\e 1\mod p$, then
$$\sum_{k=0}^{p-1}\b{2k}k\Ls u{(1-u)^2}^kf_k\e
\sum_{n=0}^{p-1}A_nu^n\mod p.$$
\par $(\t{\rm ii})$ If $u\not\e -1\mod p$, then
$$\sum_{k=0}^{p-1}\b{2k}k\Ls u{(1+u)^2}^ka_k\e
\sum_{n=0}^{p-1}A_nu^n\mod p.$$
\endpro
Proof. Taking $c_k=f_k$ in Lemma 2.1 and then applying Lemma 2.3 we
obtain (i). Taking $c_k=(-1)^ka_k$ in Lemma 2.1 and then applying
Lemma 2.3 we see that for $u\not\e 1\mod p$,
$$\sum_{k=0}^{p-1}\b{2k}k\Ls u{(1-u)^2}^k(-1)^ka_k
\e\sum_{n=0}^{p-1}u^n\cdot (-1)^nA_n\mod p.$$ Now substituting $u$
with $-u$ in the above we deduce (ii), which completes the proof.

\pro{Theorem 2.7} Let $p>3$ be a prime, $x\in\Bbb Z_p$, $x\not\e
0,-1,-\f 13\mod p$ and $\ls {9x^2+14x+9}p=1$. Then
$$\sum_{k=0}^{p-1}\b{2k}k\Ls x{9x^2+14x+9}^kf_k\e
\sum_{k=0}^{p-1}\b{2k}k^2\b{4k}{2k}\Ls x{9(1+3x)^4}^k\mod p.$$
\endpro
Proof. Let $v\in\{1,2,\ldots,\f{p-1}2\}$ be given by $v^2\e
9x^2+14x+9\mod p$, and let
$$u=\f{2x+v^2+3(x+1)v}{2x}.$$ Then $u\in\Bbb Z_p$.
Since $v^2\e 9x^2+14x+9\not\e 9(x+1)^2\mod p$ we have $v\not\e \pm
3(x+1)\mod p$. Thus $u\not\e 1\mod p$. If $u\e -1\mod p$, then
$v^2+3(x+1)v\e -4x\mod p$ and so $9(x+1)^2\e v^2+4x\e -3(x+1)v\mod
p$. As $x+1\not\e 0\mod p$ we have $v\e -3(x+1)\mod p$. We get a
contradiction. Thus $u\not\e -1\mod p$. Note that
$$\align &\f{2x+v^2+3(x+1)v}{2x}\cdot \f{2x+v^2-3(x+1)v}{2x}
\\&=\f{(2x+v^2)^2-9(x+1)^2v^2}{4x^2} \e \f{(9x^2+16x+9)^2-9(x+1)^2
(9x^2+14x+9)}{4x^2}
\\&=\f{(9x^2+16x+9)^2-(9x^2+16x+9+2x)(9x^2+16x+9-2x)}{4x^2}=1\mod p.
\endalign$$
We see that $u\not\e 0\mod p$ and
$$u+\f 1u\e \f{2x+v^2+3(x+1)v}{2x}+\f{2x+v^2-3(x+1)v}{2x}
=\f{2x+v^2}{x}\e \f{9x^2+16x+9}x\mod p.$$ Now, from the above and
Theorem 2.6 we deduce that
$$\align &\sum_{k=0}^{p-1}\b{2k}k\Ls x{9x^2+14x+9}^kf_k
\\&\e \sum_{k=0}^{p-1}\b{2k}k\f{f_k}{(u+\f 1u-2)^k}
=\sum_{k=0}^{p-1}\b{2k}k\Ls u{(1-u)^2}^kf_k
\\&\e \sum_{n=0}^{p-1}A_nu^n\e \sum_{k=0}^{p-1}\b{2k}k\Ls
u{(1+u)^2}^ka_k
\\&\e\sum_{k=0}^{p-1}\b{2k}k\f{a_k}{(u+\f 1u+2)^k}
=\sum_{k=0}^{p-1}\b{2k}k\Ls x{9(x+1)^2}^ka_k\mod p.
\endalign$$
Taking $u=\f x9$ in [S8, Theorem 4.1] we see that
$$\sum_{k=0}^{p-1}\b{2k}k\Ls x{9(x+1)^2}^ka_k
\e \sum_{k=0}^{p-1}\b{2k}k^2\b{4k}{2k}\Ls x{9(1+3x)^4}^k\mod p.$$
Thus the result follows.

\pro{Theorem 2.8} Let $p$ be a prime of the form $12k+1$ and so
$p=x^2+9y^2$. Then
$$\sum_{k=0}^{p-1}\b{2k}k\f{f_k}{(-16)^k}\e 4x^2\mod p.$$
\endpro
Proof. Taking $x=-3$ in Theorem 2.7 we see that
$$\sum_{k=0}^{p-1}\b{2k}k\f{f_k}{(-16)^k}\e \sum_{k=0}^{p-1}
\f{\b{2k}k^2\b{4k}{2k}}{(-12288)^k}\mod p.$$ Now applying [S3,
Theorem 5.3] we deduce the result.

\pro{Theorem 2.9} Let $p>5$ be a prime such that $p\e
1,5,19,23\mod{24}$. Then
$$\sum_{k=0}^{p-1}\b{2k}k\f{f_k}{96^k}
\e\cases 4x^2\mod p&\t{if $p\e 1,19\mod{24}$ and so $p=x^2+2y^2$,}
\\0\mod p&\t{if $p\e 5,23\mod{24}$.}\endcases$$
\endpro
Proof. Since $\sls 6p=1$, taking $x=9$ in Theorem 2.7 we see that
$$\sum_{k=0}^{p-1}\b{2k}k\f{f_k}{96^k}\e \sum_{k=0}^{p-1}
\f{\b{2k}k^2\b{4k}{2k}}{28^{4k}}\mod p.$$
 Now applying [S8, Theorem
5.6] we deduce the result.

\pro{Theorem 2.10} Let $p$ be a prime such that $p\e 1,9\mod{20}$
and so $p=x^2+5y^2$. Then
$$\sum_{k=0}^{p-1}\b{2k}k\f{f_k}{16^k}\e 4x^2\mod p.$$
\endpro
Proof. Let $t\in\{1,2,\ldots,\f{p-1}2\}$ be given by $t^2\e -5\mod
p$, and $x=\f{1+4t}9$. Then $x\not\e 0,-1,-\f 13\mod p$, $\f 1x\e
\f{1-4t}9\mod p$ and so $\f{9x^2+14x+9}x=14+9(x+\f 1x)=16$. Thus,
$\sls{9x^2+14x+9}p=\sls{16x}p=\sls{1+4t}p=\sls{-1-4t}p
=\sls{(2-t)^2}p=1$. We also have $\f{9(1+3x)^4}x\e -1024\mod p$.
Thus applying Theorem 2.7 and [S3, Theorem 5.5] we deduce that
$$\sum_{k=0}^{p-1}\b{2k}k\f{f_k}{16^k}\e \sum_{k=0}^{p-1}
\f{\b{2k}k^2\b{4k}{2k}}{(-1024)^k}\e 4x^2\mod p.$$ This proves the
theorem.

\pro{Theorem 2.11} Let $p>3$ be a prime and $z\in\Bbb Z_p$ with
$z\not\e \f 14\mod p$. Then
$$\sum_{n=0}^{p-1}f_nz^n\e \sum_{k=0}^{p-1}\b{2k}k\b{3k}k\Ls
z{(1-4z)^3}^k \mod p.$$
\endpro
Proof. From [Su4, (2.3)] we know that
$$f_n=\sum_{k=0}^n\b{2k}k\b{3k}k\b{n+2k}{3k}(-4)^{n-k}.$$
Thus,
$$\align \sum_{n=0}^{p-1}f_nz^n&=\sum_{n=0}^{p-1}\sum_{k=0}^n
\b{2k}k\b{3k}k\b{n+2k}{3k}(-4)^{n-k}z^n
\\&=\sum_{k=0}^{p-1}\b{2k}k\b{3k}kz^k\sum_{n=k}^{p-1}\b{n+2k}{3k}
(-4z)^{n-k}
\\&=\sum_{k=0}^{p-1}\b{2k}k\b{3k}kz^k\sum_{r=0}^{p-1-k}
\b{3k+r}{3k}(-4z)^r
\\&=\sum_{k=0}^{p-1}\b{2k}k\b{3k}kz^k
\sum_{r=0}^{p-1-k}\b{-3k-1}r(4z)^r
\\&\e \sum_{k=0}^{p-1}\b{2k}k\b{3k}kz^k
\sum_{r=0}^{p-1-k}\b{p-1-3k}r(4z)^r
\\&=\sum_{k=0}^{p-1}\b{2k}k\b{3k}kz^k(1-4z)^{p-1-3k}
\\&\e  \sum_{k=0}^{p-1}\b{2k}k\b{3k}k\Ls
z{(1-4z)^3}^k \mod p.
\endalign$$
This proves the theorem.
\par Similarly, from the formula (see [Su4, (2.2)])
$$f_n=\sum_{k=0}^{[n/2]}\b{n+k}{3k}\b{2k}k\b{3k}k2^{n-2k}$$
we deduce the following result.
 \pro{Theorem 2.12} Let $p>3$ be a prime and $z\in\Bbb
Z_p$ with $z\not\e \f 12\mod p$. Then
$$\sum_{n=0}^{p-1}f_nz^n\e \sum_{k=0}^{p-1}\b{2k}k\b{3k}k\Ls
{z^2}{(1-2z)^3}^k \mod p.$$
\endpro

\par Taking $c_k=f_k,(-1)^ka_k,\f{f_k}{(-8)^k}$ in Lemma 2.2 and
then applying Lemma 2.3 we see that for any prime $p>3$,
$$\align &\sum_{n=0}^{p-1}A_n\e \sum_{k=0}^{p-1}\f
p{2k+1}(-1)^kf_k\mod{p^3},\tag 2.3
\\&\sum_{n=0}^{p-1}(-1)^nA_n\e \sum_{k=0}^{p-1}
\f p{2k+1}a_k\mod{p^3},\tag 2.4
\\&\sum_{n=0}^{p-1}\f{D_n}{8^n}\e\sum_{k=0}^{p-1}
\f p{2k+1}\cdot \f{f_k}{8^k}\mod{p^3}.\tag 2.5
\endalign$$
It is known that ([JV])
$$f_k\e (-8)^kf_{p-1-k}\mod p\qtq{for}k=0,1,\ldots,p-1.\tag 2.6$$
Thus, from (2.3) and (2.5) we see that
$$\aligned \sum_{n=0}^{p-1}\f{D_n}{8^n}
&\e f_{\f{p-1}2}8^{-\f{p-1}2}+p\sum\Sb
k=0\\k\not=(p-1)/2\endSb^{p-1}\f 1{2k+1}\cdot\f{f_k}{8^k}
\\&=8^{-\f{p-1}2}f_{\f{p-1}2}+p\sum\Sb
k=0\\k\not=(p-1)/2\endSb^{p-1}\f 1{2(p-1-k)+1}\cdot \f{f_{p-1-k}}
{8^{p-1-k}}
\\&\e 8^{-\f{p-1}2}f_{\f{p-1}2}+p\sum\Sb
k=0\\k\not=(p-1)/2\endSb^{p-1}\f 1{-2k-1}\cdot \f{f_k/(-8)^k}
{8^{p-1-k}}
\\&\e 8^{-\f{p-1}2}f_{\f{p-1}2}-p\sum\Sb
k=0\\k\not=(p-1)/2\endSb^{p-1}\f 1{2k+1}(-1)^kf_k
\\&=((-1)^{\f{p-1}2}+8^{-\f{p-1}2})f_{\f{p-1}2}
-\sum_{k=0}^{p-1}\f{p}{2k+1}(-1)^kf_k
\\&=\f{1+(-8)^{\f{p-1}2}}{8^{\f{p-1}2}}f_{\f{p-1}2}-
\sum_{n=0}^{p-1}A_n\mod{p^2}.
\endaligned$$
If $p\e 5,7\mod 8$, then $(-8)^{(p-1)/2}\e -1\mod p$ and so $p\mid
f_{\f{p-1}2}$ by (2.6). Hence $(1+(-8)^{\f{p-1}2})f_{\f{p-1}2}\e
0\mod {p^2}$ and so $$\sum_{n=0}^{p-1}\f{D_n}{8^n}\e
-\sum_{n=0}^{p-1}A_n\mod{p^2}\qtq{for}p\e 5,7\mod 8.\tag 2.7$$

\pro{Theorem 2.13} Let $p$ be a prime with $p\e 5\mod 6$. Then
$$\sum_{n=0}^{p-1}\f{D_n}{4^n}\e 0\mod{p^2}.$$
\endpro
Proof. By [S8, Lemma 3.1],
$$D_n=\sum_{k=0}^{[n/2]}\b{2k}k^2
\b{3k}k\b{n+k}{3k}4^{n-2k}. $$ Thus,
$$\align\sum_{n=0}^{p-1}\f{D_n}{4^n}
&=\sum_{n=0}^{p-1}\sum_{k=0}^{[n/2]}\b{2k}k^2
\b{3k}k\b{n+k}{3k}4^{-2k}
\\&=\sum_{k=0}^{(p-1)/2}\f{\b{2k}k^2\b{3k}k}{16^k}\sum_{n=2k}^{p-1}
\b{n+k}{3k}
\\&=\sum_{k=0}^{(p-1)/2}\f{\b{2k}k^2\b{3k}k}{16^k}
\sum_{r=0}^{p-1-2k}\b{3k+r}r
\\&=\sum_{k=0}^{(p-1)/2}\f{\b{2k}k^2\b{3k}k}{16^k}
\sum_{r=0}^{p-1-2k}\b{-3k-1}r(-1)^r
\\&=\sum_{k=0}^{(p-1)/2}\f{\b{2k}k^2\b{3k}k}{16^k}\b{p+k}{p-1-2k}
=\sum_{k=0}^{(p-1)/2}\f{\b{2k}k^2\b{3k}k}{16^k}\b{p+k}{3k+1}.
\endalign$$
For $1\le k\le\f{p-1}2$ we see that
$$\align\b{3k}k\b{p+k}{3k+1}
&=\f p{3k+1}\cdot\f{(p^2-1^2)\cdots(p^2-k^2)(p-(k+1))\cdots(p-2k)}
{k!(2k)!}
\\&\e \f p{3k+1}\Big(1-p\sum_{r=k+1}^{2k}\f 1r\Big)\e \f p{3k+1}
\mod{p^2}.\endalign$$ Hence
$$\sum_{n=0}^{p-1}\f{D_n}{4^n}\e
\sum_{k=0}^{(p-1)/2}\f{\b{2k}k^2}{16^k} \cdot\f p{3k+1}\mod{p^2}.$$
\par For $a\in\Bbb Z_p$ let
$$S_p(a)=\sum_{k=0}^{(p-1)/2}\b ak\b{-1-a}k\f 1{3k+1}.$$ By [S9, (3.2)],
$$(3a+1)S_p(a)-(3a-1)S_p(a-1)=2\b{a-1}{\f{p-1}2}\b{-a-1}{\f{p-1}2}.$$
 For $a\not\e 0\mod
p$ we see that
$$\b{a-1}{\f{p-1}2}\b{-a-1}{\f{p-1}2}
=\f{(1^2-a^2)\cdots(\sls{p-1}2^2-a^2)}{\f{p-1}2!^2}\e 0\mod p.$$
Since $p\e 2\mod 3$, we have $p\nmid 3k+1$ for $1\le k\le \f{p-1}2$.
Hence $S_p(a)\in\Bbb Z_p$ and so
$$S_p(a)\e \f{3a-1}{3a+1}S_p(a-1)=\f{2-6a}{-2-6a}S_p(a-1)\mod p\qtq{for}a\not\e 0,-\f 13\mod
p.$$ Therefore,
$$\align S_p\Big(-\f 12\Big)&\e \f 51S_p\Big(-\f 12-1\Big)
\e \f 51\cdot\f {11}7S_p\Big(-\f 12-2\Big)
\\&\e \cdots\e \f{5\cdot 11\cdots p}{1\cdot 7\cdots(p-4)}
S_p\Big(-\f 12-\f{p+1}6\Big)\e 0\mod p.\endalign$$ Note that $\b{-\f
12}k=\b{2k}k4^{-k}$. For $p\e 2\mod 3$ we see that
$$\sum_{k=0}^{(p-1)/2}\f{\b{2k}k^2}{16^k(3k+1)}=\sum_{k=0}^{(p-1)/2}
\b{-1/2}k^2\f 1{3k+1}=S_p\Big(-\f 12\Big)\e 0\mod p$$ and so
$\sum_{n=0}^{p-1}\f {D_n}{4^n}\e 0\mod {p^2}.$ This proves the
theorem.

\pro{Theorem 2.14} Let p be an odd prime. Then
$$\sum_{k=0}^{p-1}\b{2k}k\f{a_k}{4^k}
\e\cases 4x^2\mod p&\t{if $p\e 1,3\mod 8$ and so $p=x^2+2y^2$,}
\\0\mod p&\t{if $p\e 5,7\mod 8$.}
\endcases$$
\endpro
Proof. Taking $u=1$ in Theorem 2.6(ii) we see that
$$\sum_{k=0}^{p-1}\b{2k}k\f{a_k}{4^k}\e
\sum_{n=0}^{p-1}A_n\mod p.$$ By [Su1, Corollary 1.2],
$$\sum_{n=0}^{p-1}A_n\e
\cases 4x^2\mod p&\t{if $p\e 1,3\mod 8$ and so $p=x^2+2y^2$,}
\\0\mod p&\t{if $p\e 5,7\mod 8$.}
\endcases$$
Thus the theorem is proved.
\par{\bf Remark 2.2} Let $p$ be an odd prime, and $m\in\Bbb Z_p$
with $m\not\e 0\mod p$. For conjectures on $\sum_{k=0}^{p-1}\b{2k}k
\f{a_k}{m^k}\mod {p^2}$, see [Su3, Conjectures 7.8 and 7.9] and [S8,
Conjectures 6.4-6.6].

\pro{Conjecture 2.1} Let $p$ be an odd prime. If $p\e 1,3\mod 8$ and
so $p=x^2+2y^2$, then
$$f_{\f{p-1}2}\e (-1)^{\f{p-1}2}\big((3\cdot
2^{p-1}+1)x^2-2p)\mod{p^2}$$  and
$$f_{\f{p^2-1}2}\e 4x^4(3\cdot 2^{p-1}+1)-16x^2p\mod{p^2}.$$
\endpro

\pro{Conjecture 2.2} Let $p$ be an odd prime. If $p\e 5,7\mod 8$,
then
$$f_{\f{p^2-1}2}\e p^2\mod{p^3}\qtq{and}f_{\f{p^r-1}2}
\e 0\mod{p^r}\qtq{for} r\in\Bbb Z^+.$$
\endpro

\section*{3. Congruences involving $\{S_n\}$}
\par Recall that
$$S_n(x)=\sum_{k=0}^n\b nk\b{2k}k\b{2n-2k}{n-k}x^k
\q (n=0,1,2,\ldots)$$ and $S_n=S_n(1)$. From [G, (6.12)] we know
that
$$S_n(-1)=\sum_{k=0}^n\b nk\b{2k}k\b{2n-2k}{n-k}(-1)^k=\cases 0&\t{if $n$ is odd,}
\\\b n{n/2}^2&\t{if $n$ is even.}
\endcases\tag 3.1$$
Using Maple and  the Zeilberger algorithm we see that
$$n^2S_n=4(3n^2-3n+1)S_{n-1}-32(n-1)^2S_{n-2}\q (n\ge 2).$$

\pro{Lemma 3.1} For any nonnegative integer $n$ we have
$$S_n(-x)=\sum_{k=0}^n\b nk(-1)^k4^{n-k}S_k(x).$$
\endpro
Proof. Since $\b{-\f 12}k=\f{\b{2k}k}{(-4)^k}$, using Vandermonde's
identity we see that for any nonnegative integer $m$,
$$\align \sum_{k=0}^m\b mk\b{2k}k(-1)^k4^{m-k}
&=4^m\sum_{k=0}^m\b m{m-k}\b{-\f 12}k =4^m\b{m-\f 12}m\\&=4^m\cdot
(-1)^m\b{-\f 12}m=\b{2m}m.\endalign$$ Note that $\b nk \b kr=\b
nr\b{n-r}{k-r}$. From the above we see that
$$\align\sum_{k=0}^n\b nk(-1)^k4^{n-k}S_k(x)
&=\sum_{k=0}^n\b nk(-1)^k4^{n-k}\sum_{r=0}^k\b
kr\b{2r}r\b{2(k-r)}{k-r}x^r
\\&=\sum_{r=0}^n\b{2r}rx^r\b
nr\sum_{k=r}^n\b{n-r}{k-r}\b{2(k-r)}{k-r} (-1)^k4^{n-k}
\\&=\sum_{r=0}^n\b nr\b{2r}rx^r(-1)^r\sum_{s=0}^{n-r}\b{n-r}s\b {2s}s(-1)^s4^{n-r-s}
\\&=\sum_{r=0}^n\b nr\b{2r}r(-x)^r\cdot \b{2n-2r}{n-r}=S_n(-x).
\endalign$$
This proves the lemma.

\pro{Lemma 3.2} For any nonnegative integer $m$ we have
$$\sum_{k=0}^m\b mkS_k(x)n^{m-k}
=\sum_{k=0}^m\b mk(-1)^kS_k(-x)(n+4)^{m-k}$$ and so
$$\sum_{k=0}^m \b mkS_kn^{m-k}=\sum_{k=0}^{[m/2]}\b
m{2k}\b{2k}k^2(n+4)^{m-2k}.$$

\endpro
Proof. Note that $\b mk\b kr=\b mr\b{m-r}{k-r}$. By Lemma 3.1,
$$\align \sum_{k=0}^m\b mkS_k(x)n^{m-k}&
=\sum_{k=0}^m\b mkn^{m-k}\sum_{r=0}^k\b kr(-1)^rS_r(-x)4^{k-r}
\\&=\sum_{r=0}^m(-1)^rS_r(-x)\sum_{k=r}^m\b mk\b kr4^{k-r}n^{m-k}
\\&=\sum_{r=0}^m\b mr(-1)^rS_r(-x)n^{m-r}\sum_{k=r}^m\b{m-r}{k-r}
\Ls 4n^{k-r} \\&=\sum_{r=0}^m\b mr(-1)^rS_r(-x)n^{m-r}\Big(1+\f
4n\Big)^{m-r}
\\&=\sum_{r=0}^m\b mr(-1)^rS_r(-x)(n+4)^{m-r}.\endalign$$
Taking $x=1$ in the above formula and then applying (3.1) we deduce
the remaining result.
\par If $\{c_n\}$ is a sequence satisfying
$$\sum_{k=0}^n\binom nk(-1)^kc_k=c_n\q (n=0,1,2,\ldots),$$
we say that $\{c_n\}$ is an even sequence. In [S1,S6] the author
investigated the properties of even sequences. \pro{Lemma 3.3}
Suppose that $\{c_n\}$ is an even sequence.
\par $(\t{\rm i})$ $(\t{\rm [S7, Theorem 2.3])}$ If $n$ is odd, then
$$\sum_{k=0}^n\b nk\b{n+k}k(-1)^kc_k=0.$$
\par $(\t{\rm ii})$ $(\t{\rm [S7, Theorems 5.3 and 5.4]})$ If $p$ is a prime of the form $4k+3$
and $c_0,c_1,\ldots,c_{\f{p-1}2}\in\Bbb Z_p$, then
$$\sum_{k=0}^{(p-1)/2}\b{2k}k^2\f{c_k}{16^k}
\e 0\mod{p^2}\qtq{and}\sum_{k=0}^{(p-1)/2}\b{2k}k\f{c_k}{2^k}\e
0\mod p.$$
\endpro

 \pro{Theorem  3.1} Let $n$ be a nonnegative integer. Then
$$\align
&\sum_{k=0}^n\b nk(-1)^k4^{n-k}S_k=\cases 0&\t{if $n$ is odd,}
\\\b n{n/2}^2&\t{if $n$ is even,}\endcases\tag i
\\&\sum_{k=0}^n\b nk(-1)^k\f{S_k}{8^k}=\f{S_n}{8^n},\tag ii
\\&\sum_{k=0}^n\b nk\b {n+k}k(-8)^{n-k}S_k= \cases 0&\t{if $n$ is
odd,}
\\(-1)^{\f n2}\b n{n/2}^3&\t{if $n$ is even.}
\endcases\tag iii\endalign $$
\endpro
Proof. Taking $x=1$ in Lemma 3.1 and then applying (3.1) we deduce
part (i). By Lemma 3.2,
$$\align &\sum_{k=0}^n\b nkS_km^{n-k}
\\&=\sum_{k=0}^{[n/2]}\b n{2k}\b{2k}k^2(m+4)^{n-2k}
=(-1)^n\sum_{k=0}^{[n/2]}\b n{2k}\b{2k}k^2(-m-4)^{n-2k}
\\&=(-1)^n\sum_{k=0}^n\b nkS_k(-m-8)^{n-k}
.\endalign$$ That is, $$\sum_{k=0}^n\b nkS_km^{n-k}=\sum_{k=0}^n\b
nk(-1)^kS_k(m+8)^{n-k}.\tag 3.2$$ Putting $m=0$ in (3.2) we obtain
part (ii). By (ii), $\{\f{S_n}{8^n}\}$ is an even sequence. Thus
applying Lemma 3.3(i), for odd $n$ we have
 $$\sum_{k=0}^n\b nk\b
{n+k}k(-8)^{n-k}S_k=(-8)^n\sum_{k=0}^n\b nk\b
{n+k}k(-1)^k\f{S_k}{8^k}=0.$$
 Let
$$c_n=\sum_{k=0}^n\b nk\b {n+k}k(-8)^{n-k}S_k.$$ Then $c_0=1$. Using Maple software doublesum.mpl
and the method in [CHM] we find that $c_n=(\f {4(n-1)}n)^3c_{n-2}$.
 When $n$ is even we see that
$$\f{(-1)^{n/2}\b n{n/2}^3}{(-1)^{(n-2)/2}\b{n-2}{(n-2)/2}^3}
=\Ls{4(n-1)}n^3.$$ Thus part (iii) holds for even $n$. The proof is
now complete.

\pro{Lemma 3.4} Let $p$ be an odd prime, $x\in\Bbb Z_p$ and $x\not\e
-1\mod p$. Then
$$\sum_{k=0}^{p-1}\b{2k}k\Ls x{8(1+x)^2}^kS_k\e
\sum_{k=0}^{\f{p-1}2}\b{2k}k^3\Big(-\f{x^2}{64}\Big)^k\mod p.$$
\endpro
Proof. Taking $u=-x$ and $c_k=\f{S_k}{(-8)^k}$ in Lemma 2.1 and then
applying Theorem 3.1(iii) we see that
$$\align &\sum_{k=0}^{p-1}\b{2k}k\Ls{-x}{(1+x)^2}\f{S_k}{(-8)^k}
\\&\e \sum_{n=0}^{p-1}(-x)^n\sum_{k=0}^n\b nk\b{n+k}k\f{S_k}{(-8)^k}
=\sum_{k=0}^{(p-1)/2}(-x)^{2k}\cdot \f{(-1)^k}{(-8)^{2k}}\b{2k}k^3
\mod p.\endalign$$ This yields the result.

\pro{Theorem 3.2} Let $p$ be an odd prime. Then
$$\sum_{k=0}^{p-1}\b{2k}k\f{S_k}{32^k}
\e\cases 4x^2\mod p&\t{if $p\e 1,3\mod 8$ and so $p=x^2+2y^2$,}
\\0\mod p&\t{if $p\e 5,7\mod 8$.}
\endcases$$
\endpro
Proof. Taking $x=1$ in Lemma 3.4 we find that
$$\sum_{k=0}^{p-1}\b{2k}k\f{S_k}{32^k}\e
\sum_{k=0}^{(p-1)/2}\b{2k}k^3\f 1{(-64)^k}\mod p.$$ Now applying
[S3, Theorems 3.3-3.4] we deduce the result.

\pro{Theorem 3.3} Let $p$ be an odd prime. Then
$$\sum_{k=0}^{p-1}\b{2k}k\f{S_k}{16^k}\e
\cases 4x^2\mod p&\t{if $p\e 1\mod 4$ and so $p=x^2+4y^2$,} \\0\mod
p&\t{if $p\e 3\mod 4$.}\endcases$$
\endpro
Proof. From Theorem 3.1(ii) we know that $\{\f{S_n}{8^n}\}$ is an
even sequence. Thus applying Lemma 3.3(ii) we have
$\sum_{k=0}^{p-1}\b{2k}k\f{S_k}{16^k}\e 0\mod p$ for $p\e 3\mod 4$.
Now assume $p\e 1\mod 4$ and so $p=x^2+4y^2$. Let
$t\in\{1,2,\ldots,\f{p-1}2\}$ be given by $t^2\e -1\mod p$. By Lemma
3.4,
$$\align \sum_{k=0}^{p-1}\b{2k}k\f{S_k}{16^k}
&\e \sum_{k=0}^{p-1}\b{2k}k\Ls t{8(1+t)^2}^kS_k \e
\sum_{k=0}^{(p-1)/2}\b{2k}k^3\Big(-\f{t^2}{64}\Big)^k
\\&\e\sum_{k=0}^{(p-1)/2}\b{2k}k^3\f 1{64^k}\mod p.\endalign$$
It is well known that (see for example [Ah])
$$\sum_{k=0}^{(p-1)/2}\b{2k}k^3\f 1{64^k} \e 4x^2-2p \mod {p^2}.$$
Thus $\sum_{k=0}^{p-1}\b{2k}k\f{S_k}{16^k}\e 4x^2\mod p$, which
completes the proof.

\pro{Theorem 3.4} Let $p$ be an odd prime and
$q_p(2)=(2^{p-1}-1)/p$. Then
$$\sum_{k=0}^{p-1}\b{2k}k^2\f{S_k}{128^k}
\e\cases (-1)^{\f{p-1}4}(8x^3+6x(2q_p(2)x^2-1)p)\mod {p^2}
\\\qq\qq\qq\q\t{if $p=x^2+4y^2\e 1\mod 4$ and $4\mid x-1$,}
\\0\mod {p^2}\qq\; \t{if $p\e 3\mod 4$.}\endcases$$
\endpro
Proof. It is clear that for $k\in\{0,1,\ldots,\f{p-1}2\}$,
$$\aligned &\b{\f{p-1}2}k\b{\f{p-1}2+k}k\\&=\b{2k}k
\b{\f{p-1}2+k}{2k}=\b{2k}k
\f{(p^2-1^2)(p^2-3^2)\cdots(p^2-(2k-1)^2)}{2^{2k}\cdot (2k)!} \\& \e
\b{2k}k(-1)^k\f{ 1^2\cdot 3^2\cdots (2k-1)^2}{2^{2k}\cdot (2k)!}
=\f{\b{2k}k^2}{(-16)^k}\mod{p^2}.
\endaligned\tag 3.3$$
Thus, from Theorem 3.1(iii) we deduce that $$
\align\sum_{k=0}^{p-1}\b{2k}k^2\f{S_k}{128^k} &\e
\sum_{k=0}^{(p-1)/2}\b{\f{p-1}2}k\b{\f{p-1}2+k}k\f{S_k}{8^k}\\&
=\cases 0\mod{p^2}&\t{if $p\e 3\mod 4$,}
\\(-1)^{\f{p-1}4}\b{\f{p-1}2}{\f{p-1}4}^3\mod{p^2}
&\t{if $p\e 1\mod 4$.}\endcases\endalign$$ By [CDE], for
$p=x^2+4y^2\e 1\mod 4$ with $4\mid x-1$,
$$\align \b{\f{p-1}2}{\f{p-1}4}&\e \f{2^{p-1}+1}2\Big(2x-\f p{2x}\Big)
=\Big(1+\f 12 q_p(2)p\Big)\Big(2x-\f p{2x}\Big) \\& \e
2x+p\Big(q_p(2)x-\f 1{2x}\Big)\mod{p^2}.\endalign$$ Thus,
$$ \b{\f{p-1}2}{\f{p-1}4}^3\e
\Big(2x+p\Big(q_p(2)x-\f 1{2x}\Big)\Big)^3\e
8x^3+6x(2q_p(2)x^2-1)p\mod {p^2}.$$ Now putting the above together
we deduce the result.
\par For an odd prime $p$ and $a\in\Bbb Z_p$ let
$\ap\in\{0,1,\ldots,p-1\}$ be given by $a\e \ap\mod p$.

\pro{Theorem 3.5} Let $p>3$ be a prime, $a\in\Bbb Z_p$ and $\ap\e
1\mod 2$. Then
$$\sum_{k=0}^{p-1}\b ak\b{-1-a}k\f{S_k}{8^k}\e 0\mod{p^2}.$$
In particular, for $a=-\f 13,-\f 14,-\f 16$ we have
$$\align &\sum_{k=0}^{p-1}\f{\b{2k}k\b{3k}k}{216^k}S_k\e 0\mod{p^2}
\qtq{for}p\e 2\mod 3,
\\&\sum_{k=0}^{p-1}\f{\b{2k}k\b{4k}{2k}}{512^k}S_k\e 0\mod{p^2}
\qtq{for}p\e 3\mod 4,
\\&\sum_{k=0}^{p-1}\f{\b{3k}k\b{6k}{3k}}{3456^k}S_k\e 0\mod{p^2}
\qtq{for}p\e 2\mod 3.\endalign$$
\endpro
Proof. This is immediate from Theorem 3.1(ii) and [S5, Theorem 2.4].

 \pro{Theorem
3.6} Let $p$ be an odd prime, $n\in\Bbb Z_p$ and $n\not\e 0,-16\mod
p$. Then
$$\sum_{k=0}^{p-1}\b{2k}k\f{S_k}{(n+16)^k}\e \Ls {n(n+16)}p
\sum_{k=0}^{p-1}
\f{\b{2k}k^2\b{4k}{2k}}{n^{2k}}\mod p.$$
\endpro
Proof. Clearly $p\mid \b{2k}k$ for $\f p2<k<p$ and $p\mid
\b{2k}k\b{4k}{2k}$ for $\f p4<k<p$. Note that $\b{\f{p-1}2}k \e
\b{-\f 12}k=\f{\b{2k}k}{(-4)^k}\mod p$ for $0\le k\le\f{p-1}2$.  By
Lemma 3.2,
$$\align &\sum_{k=0}^{p-1}\b{2k}k\f{S_k}{(n+16)^k}
\\&\e \sum_{k=0}^{\f{p-1}2}\b{\f{p-1}2}kS_k\Ls{-4}{n+16}^k
\e \Ls{-n-16}p \sum_{k=0}^{\f{p-1}2}\b{\f{p-1}2}kS_k
\Ls{n+16}{-4}^{\f{p-1}2-k}
\\&= \Ls{-n-16}p \sum_{k=0}^{[p/4]}\b{\f{p-1}2}{2k}
\b{2k}k^2\Big(-\f n4\Big)^{\f{p-1}2-2k} \\&\e
\Ls{-n(-n-16)}p\sum_{k=0}^{[p/4]}
\f{\b{4k}{2k}}{(-4)^{2k}}\b{2k}k^2\f 1{(-n/4)^{2k}} \e \Ls{n(n+16)}p
\sum_{k=0}^{p-1} \f{\b{2k}k^2\b{4k}{2k}}{n^{2k}}\mod p.\endalign$$
This proves the theorem.

\pro{Theorem 3.7} Let $p>7$ be a prime. Then
$$\align &\sum_{k=0}^{p-1}\b{2k}k\f{S_k}{7^k}\e
\sum_{k=0}^{p-1}\b{2k}k\f{S_k}{25^k}
\\&\e\cases 4x^2\mod p&\t{if $p\e 1,2,4\mod 7$ and so $p=x^2+7y^2$,}
\\0\mod p&\t{if $p\e 3,5,6\mod 7$.}
\endcases\endalign$$
\endpro
Proof. Taking $n=\pm 9$ in Theorem 3.6 and then applying [S3,
Theorem 5.2] we deduce the result.

\pro{Theorem 3.8} Let $p$ be a prime such that $p\e 1,7,17,23
\mod{24}$. Then
$$\align &\Ls 3p\sum_{k=0}^{p-1}\b{2k}k\f{S_k}{64^k}\e
\Ls 6p\sum_{k=0}^{p-1}\b{2k}k\f{S_k}{(-32)^k}
\\&\e\cases 4x^2\mod p&\t{if $p\e 1,7\mod {24}$ and so $p=x^2+6y^2$,}
\\0\mod p&\t{if $p\e 17,23\mod {24}$.}
\endcases\endalign$$
\endpro
Proof. Taking $n=\pm 48$ in Theorem 3.6 and then applying [S3,
Theorem 5.4] we deduce the result.

\pro{Theorem 3.9} Let $p>5$ be a prime.  Then
$$\align &\Ls 2p\sum_{k=0}^{p-1}\b{2k}k\f{S_k}{800^k}\e
\Ls 3p\sum_{k=0}^{p-1}\b{2k}k\f{S_k}{(-768)^k}
\\&\e\cases 4x^2\mod p&\t{if $p\e 1,3\mod {8}$ and so $p=x^2+2y^2$,}
\\0\mod p&\t{if $p\e 5,7\mod {8}$.}
\endcases\endalign$$
\endpro
Proof. By [S8, Theorem 5.6],
$$\sum_{k=0}^{p-1}\f{\b{2k}k^2\b{4k}{2k}}{28^{4k}}
\e\cases 4x^2\mod p&\t{if $p=x^2+2y^2\e 1,3\mod 8$,}
 \\0\mod {p^2}&\t{if $p\e 5,7\mod 8$}
 \endcases$$
 Now taking $n=\pm 28^2=\pm 784$ in Theorem 3.6 and then applying
 the above we obtain the result.

 \pro{Theorem 3.10} Let $p$ be a prime such that $p\e 1,9\mod{10}$. Then
$$\align &\sum_{k=0}^{p-1}\b{2k}k\f{S_k}{160^k}\e
\sum_{k=0}^{p-1}\b{2k}k\f{S_k}{(-128)^k}
\\&\e\cases \sls 2p4x^2\mod p&\t{if $p\e 1,9,11,19\mod {40}$
 and so $p=x^2+10y^2$,}
\\0\mod p&\t{if $p\e 21,29,31,39\mod {40}$.}
\endcases\endalign$$
\endpro
Proof. Taking $n=\pm 144$ in Theorem 3.6 and then applying [S3,
(5.9)] we deduce the result.

\pro{Theorem 3.11} Let $p$ be a prime such that $p\e \pm 1\mod{8}$.
Then
$$\align &\sum_{k=0}^{p-1}\b{2k}k\f{S_k}{1600^k}\e
\sum_{k=0}^{p-1}\b{2k}k\f{S_k}{(-1568)^k}
\\&\e\cases \sls {-1}p4x^2\mod p&\t{if $p\e 1,3,4,5,9\mod {11}$
 and so $p=x^2+22y^2$,}
\\0\mod p&\t{if $p\e 2,6,7,8,10\mod {11}$.}
\endcases\endalign$$
\endpro
Proof. Taking $n=\pm 1584$ in Theorem 3.6 and then applying [S3,
(5.9)] we deduce the result.

\pro{Theorem 3.12} Let $p$ be a prime such that $\sls p{29}=1$. Then
$$\align &\sum_{k=0}^{p-1}\b{2k}k\f{S_k}{156832^k}\e
\sum_{k=0}^{p-1}\b{2k}k\f{S_k}{(-156800)^k}
\\&\e\cases \sls 2p4x^2\mod p&\t{if $p\e 1,3\mod {8}$
 and so $p=x^2+58y^2$,}
\\0\mod p&\t{if $p\e 5,7\mod {8}$.}
\endcases\endalign$$
\endpro
Proof. Taking $n=\pm 396^2=\pm 156816$ in Theorem 3.6 and then
applying [S3, (5.9)] we deduce the result.

\pro{Theorem 3.13} Let $p$ be an odd prime, $n\in\Bbb Z_p$ and
$n\not\e 0\mod p$.
\par $(\t{\rm i})$ If $n\not\e 4\mod p$, then
$$\sum_{k=0}^{p-1}\f{S_k(x)}{n^k}\e \sum_{k=0}^{p-1}
\f{S_k(-x)}{(4-n)^k}\mod p.$$

\par $(\t{\rm ii})$ If $n\not\e 16\mod p$, then
$$\sum_{k=0}^{p-1}\b{2k}k\f{S_k(x)}{n^k}\e \Ls{n(n-16)}p
\sum_{k=0}^{p-1}\b{2k}k \f{S_k(-x)}{(16-n)^k}\mod p.$$
\endpro
Proof. Since $\b{p-1}k\e (-1)^k\mod p$ and $\b{\f{p-1}2}k\e
\f{\b{2k}k}{(-4)^k}\mod p$, taking $m=p-1$ and replacing $n$ with
$-n$ in Lemma 3.2 we deduce part (i), and taking $m=\f{p-1}2$ and
replacing $n$ with $-\f n4$ in Lemma 3.2 we deduce part (ii).

\pro{Conjecture 3.1} Let $p$ be an odd prime, $n\in \{\pm 156816,\pm
1584,\pm 784,\pm 144,\pm 48,\pm 9\}$ and $n\not\e 0,\pm 16\mod p$.
Then
$$\sum_{k=0}^{p-1}\b{2k}k\f{S_k}{(n+16)^k}\e \Ls {n(n+16)}p\sum_{k=0}^{p-1}
\f{\b{2k}k^2\b{4k}{2k}}{n^{2k}}\mod {p^2}.$$
\endpro

\pro{Conjecture 3.2} Let $p$ be an odd prime. If $p\e 1,3\mod 8$ and
so $p=x^2+2y^2$, then
$$S_{\f{p-1}2}\e (5\cdot
2^{p-1}-1)x^2-2p\mod{p^2}$$ and
$$S_{\f{p^2-1}2}\e 4x^4(5\cdot 2^{p-1}-1)-16x^2p\mod{p^2}.$$
\endpro

\pro{Conjecture 3.3} Let $p$ be an odd prime. If $p\e 5,7\mod 8$,
then
$$S_{\f{p^2-1}2}\e p^2\mod{p^3}\qtq{and}S_{\f{p^r-1}2}
\e 0\mod{p^r}\qtq{for} r\in\Bbb Z^+.$$
\endpro

\pro{Conjecture 3.4} For $m=2,3,4,\ldots$ we have
$$S_m^2<S_{m+1}S_{m-1}<\Big(1+\f 1{m(m-1)}\Big)S_m^2.$$
\endpro

\section*{4. Congruences involving $\{C_n\}$}
\par For any nonnegative integer $n$ define
$$C_n(x)=\sum_{k=0}^{[n/3]}\b{2k}k\b{3k}k\b n{3k}x^{n-3k}
\qtq{and}C_n=C_n(-3).$$
 \pro{Lemma 4.1} Let $m$ be a nonnegative
integer. Then
$$\sum_{k=0}^m\b mkC_k(x)n^{m-k}=C_m(x+n).$$
\endpro
Proof. It is clear that
$$\align \sum_{k=0}^m\b mkC_k(x)n^{m-k}&=\sum_{k=0}^m
\b mkn^{m-k}\sum_{r=0}^k\b{2r}r\b{3r}r\b k{3r}x^{k-3r}
\\&=\sum_{r=0}^m\b{2r}r\b{3r}rn^{m-3r}
\sum_{k=r}^m\b mk\b k{3r}\Ls xn^{k-3r}
\\&=\sum_{r=0}^m\b{2r}r\b{3r}rn^{m-3r}\sum_{k=3r}^m\b
m{3r}\b{m-3r}{k-3r}\Ls xn^{k-3r}
\\&=\sum_{r=0}^m\b{2r}r\b{3r}r\b m{3r}n^{m-3r}
\Big(1+\f xn\Big)^{m-3r}=C_m(x+n).
\endalign$$ This proves the lemma.
\pro{Theorem 4.1} Let $p$ be an odd prime, $n,x\in\Bbb Z_p$ and
$n(n+4x)\not\e 0\mod p$. Then
$$\sum_{k=0}^{p-1}\b{2k}k\f{C_k(x)}{(n+4x)^k}
\e\Ls {n(n+4x)}p\sum_{k=0}^{p-1}\f{\b{2k}k\b{3k}k\b{6k}{3k}}
{n^{3k}}\mod p.$$
\endpro
Proof. As $\b{\f{p-1}2}k\e \b{2k}k4^{-k}\mod p$ and $p\mid \b{2k}k
\b{3k}k\b{6k}{3k}$ for $\f p6<k<p$, using Lemma 4.1 we see that
$$\align &\sum_{k=0}^{p-1}\b{2k}k\f{C_k(x)}{(n+4x)^k}
\\&\e \sum_{k=0}^{(p-1)/2}\b{\f{p-1}2}kC_k(x)\Ls {-4}{n+4x}^k
\e \Ls{-4(n+4x)}p\sum_{k=0}^{(p-1)/2}\b{\f{p-1}2}kC_k(x)
\Ls{n+4x}{-4}^{\f{p-1}2-k}
\\&=\Ls{-n-4x}pC_{\f{p-1}2}\Big(-\f n4\Big)
=\Ls{-n-4x}p\sum_{k=0}^{[p/6]}\b{2k}k\b{3k}k\b{\f{p-1}2}{3k}
\Big(-\f n4\Big)^{\f{p-1}2-3k}
\\&\e \Ls{n(n+4x)}p\sum_{k=0}^{[p/6]}\b{2k}k\b{3k}k\b{6k}{3k}\f
1{(-4)^{3k}\cdot (-n/4)^{3k}}
\\&\e \Ls {n(n+4x)}p\sum_{k=0}^{p-1}\f{\b{2k}k\b{3k}k\b{6k}{3k}}
{n^{3k}}\mod p.
\endalign$$
This proves the theorem.

\pro{Theorem 4.2} Let $p$ be a prime, $p\not=2,3,11$, $t\in\Bbb Z_p$
and $33+2t\not\e 0\mod p$. Then
$$\Ls {33+2t}p\sum_{k=0}^{p-1}\b{2k}k\f{C_k(t)}{(66+4t)^k}
\e\cases 4x^2\mod p&\t{if $p=x^2+4y^2\e 1\mod 4$,}
\\0\mod p&\t{if $p\e 3\mod 4$.}
\endcases$$\endpro
Proof. Taking $n=66$ and replacing $x$ with $t$ in Theorem 4.1 and
then applying [S4, Theorem 4.3] we deduce the result.

\pro{Theorem 4.3} Let $p>5$ be a prime, $t\in\Bbb Z_p$ and $t\not\e
-5\mod p$. Then
$$\Ls {-(5+t)}p\sum_{k=0}^{p-1}\b{2k}k\f{C_k(t)}{(20+4t)^k}
\e\cases 4x^2\mod p&\t{if $p=x^2+2y^2\e 1,3\mod 8$,}
\\0\mod p&\t{if $p\e 5,7\mod 8$.}
\endcases$$\endpro
Proof. Taking $n=20$ and replacing $x$ with $t$ in Theorem 4.1 and
then applying [S4, Theorem 4.4] we deduce the result.

\pro{Theorem 4.4} Let $p>7$ be a prime, $t\in\Bbb Z_p$ and $4t\not\e
15\mod p$. Then
$$\align \Ls {-15+4t}p\sum_{k=0}^{p-1}\b{2k}k\f{C_k(t)}{(-15+4t)^k}
\e\cases 4x^2\mod p&\t{if $p=x^2+7y^2\e 1,2,4\mod 7$,}
\\0\mod p&\t{if $p\e 3,5,6\mod 7$.}
\endcases\endalign$$\endpro
Proof. Taking $n=-15$ and replacing $x$ with $t$ in Theorem 4.1 and
then applying [S4, Theorem 4.7] we deduce the result.

\pro{Theorem 4.5} Let $p>7$ be a prime, $t\in\Bbb Z_p$ and $4t\not\e
-255\mod p$. Then
$$\align \Ls {-255-4t}p\sum_{k=0}^{p-1}\b{2k}k\f{C_k(t)}{(255+4t)^k}
\e\cases 4x^2\mod p&\t{if $p=x^2+7y^2\e 1,2,4\mod 7$,}
\\0\mod p&\t{if $p\e 3,5,6\mod 7$.}
\endcases\endalign$$\endpro
Proof. Taking $n=255$ and replacing $x$ with $t$ in Theorem 4.1 and
then applying [S4, Theorem 4.7] we deduce the result.

\pro{Theorem 4.6} Let $p$ be a prime,  $p\not=2,3,11$, $t\in\Bbb
Z_p$ and $t\not\e 8\mod p$. Then
$$\align \Ls {t-8}p\sum_{k=0}^{p-1}\b{2k}k\f{C_k(t)}{(4t-32)^k}
\e\cases x^2\mod p&\t{if $\sls p{11}=1$ and so $4p=x^2+11y^2$,}
\\0\mod p&\t{if $\sls p{11}=-1$.}
\endcases\endalign$$\endpro
Proof. Taking $n=-32$ and replacing $x$ with $t$ in Theorem 4.1 and
then applying [S4, Theorem 4.8] we deduce the result.

\pro{Theorem 4.7} Let $p$ be a prime,  $p\not=2,3,19$, $t\in\Bbb
Z_p$ and $t\not\e 24\mod p$. Then
$$\align \Ls {t-24}p\sum_{k=0}^{p-1}\b{2k}k\f{C_k(t)}{(4t-96)^k}
\e\cases x^2\mod p&\t{if $\sls p{19}=1$ and so
 $4p=x^2+19y^2$,}
\\0\mod p&\t{if $\sls p{19}=-1$.}
\endcases\endalign$$\endpro
Proof. Taking $n=-96$ and replacing $x$ with $t$ in Theorem 4.1 and
then applying [S4, Theorem 4.9] we deduce the result.
\par Using Theorem 4.1 and [S4, Theorem 4.9] one can also deduce the
following results.

\pro{Theorem 4.8} Let $p$ be a prime, $p\not=2,3,5,43$, $t\in\Bbb
Z_p$ and $t\not\e 240\mod p$. Then
$$\align \Ls {t-240}p\sum_{k=0}^{p-1}\b{2k}k\f{C_k(t)}{(4t-960)^k}
\e\cases x^2\mod p&\t{if $\sls p{43}=1$ and so
 $4p=x^2+43y^2$,}
\\0\mod p&\t{if $\sls p{43}=-1$.}
\endcases\endalign$$\endpro

\pro{Theorem 4.9} Let $p$ be a prime, $p\not=2,3,5,11,67$, $t\in\Bbb
Z_p$ and $t\not\e 1320\mod p$. Then
$$\align \Ls {t-1320}p\sum_{k=0}^{p-1}\b{2k}k\f{C_k(t)}{(4t-5280)^k}
\e\cases x^2\mod p&\t{if $\sls p{67}=1$ and so
 $4p=x^2+67y^2$,}
\\0\mod p&\t{if $\sls p{67}=-1$.}
\endcases\endalign$$\endpro

\pro{Theorem 4.10} Let $p$ be a prime, $p\not=2,3,5,23,29,163$,
$t\in\Bbb Z_p$ and $t\not\e 160080\mod p$. Then
$$\align &\Ls {t-160080}p\sum_{k=0}^{p-1}\b{2k}k
\f{C_k(t)}{(4t-640320)^k}
\\&\e\cases x^2\mod p&\t{if $\sls p{163}=1$ and so
 $4p=x^2+163y^2$,}
\\0\mod p&\t{if $\sls p{163}=-1$.}
\endcases\endalign$$\endpro

\pro{Conjecture 4.1} Let $p$ be an odd prime,
 $n\in\{-640320,-5280,-960,-96,-32,-15,$ $20,66,255\}$ and
$n(n-12)\not\e 0\mod p$. Then
$$\sum_{k=0}^{p-1}\b{2k}k\f{C_k}{(n-12)^k}
\e\Ls {n(n-12)}p\sum_{k=0}^{p-1}\f{\b{2k}k\b{3k}k\b{6k}{3k}}
{n^{3k}}\mod {p^2}.$$
\endpro

\end{document}